\documentclass[12pt]{article}
\usepackage{amsmath}
\usepackage[T2A]{fontenc}
\usepackage[cp1251]{inputenc}
\usepackage[english, russian]{babel}
\usepackage[12pt]{extsizes}
\usepackage[mag=1000,letterpaper,
left=2cm,right=2cm,top=3cm,
bottom=4cm]{geometry}
\usepackage[dvips]{graphicx}
\begin{document}
%\tableofcontents\newpage
%Mathematical Notes, Vol. 61,
%No. 5, 1998

\vskip 14pt\noindent

\begin{center}{\bf On the 
Asymptotic Integration of a 
System of Linear Differential 
Equations with Oscillatory 
Decreasing Coefficients}
\end{center}

\vskip 14pt\noindent

\centerline{V.Sh. Burd and 
V.A. Karakulin} 

\vskip 18pt\noindent

Abstact. A system of linear differential 
equations with oscillatory 
decreasing coefficients is considered.
The coefficients has the form
$t^{-\alpha}a(t)$,~$\alpha>0$,
where $a(t)$ is trigonometric
polynomial with an arbitrary set
of frequencies. The asymptotic
behavior of the solutions of 
this system as $t\to\infty$ is
studied. We construct an invertible
(for sufficiently large $t$)
change of variables that takes
the original system to a system
not containing oscillatory
coefficients in its principal
part. The study of the asymptotic
behavior of the solutions of the
transformed system is a simpler
problem. As an example, the
following equation is considered:
$$
\frac{d^2x}{dt^2}+\left(1+\frac{\sin\lambda t}
{t^\alpha}\right)x=0,
$$
where $\lambda$ and $\alpha$,~
$0<\alpha\le 1$, are real numbers.

\vskip 14pt\noindent

The stability and asymptotic 
behavior of the solutions of 
linear systems
of differential equations
$$
\frac{dx}{dt}=Ax+B(t)x,
$$
where $A$ is a constant matrix, and matrix 
$B(t)$ is small in a certain 
sense when $t\to\infty$ has 
been studied
by many authors
(see [1-6]) as well as the 
monographs
[7-11]).
\par
Shtokalo [12-13] studied
of the stability of
solutions of the following 
system of differential equations
$$
\frac{dx}{dt}=Ax+
\varepsilon B(t)x.
\eqno(1)  
$$
Here $\varepsilon>0$ is a small 
parameter, 
$A$ is a constant square matrix, 
and
$B(t)$ is a square matrix whose 
elements are
trigonometric polynomials 
$b_{kl}(t)$~$(k,l =1,...,m)$ in the 
form
$$
b_{kl}(t)=\sum_{j=1}^m b_{j}^{kl}
e^{i\lambda_j t};  
$$
Matrices with such elements are
called matrices of class 
$Sigma$. 
The mean of a matrix from  $\Sigma$ 
is a constant matrix
that consists of the constant 
terms 
$(\lambda_j=0)$ of the elements
of the matrix.
\par
Using the Bogoliubov averaged
method, Shtokalo transformed 
system (1) into a system with
constant coefficients depending 
on the parameter $\varepsilon$
up to terms of any order in 
smallnes in $\varepsilon$.
\par
We also mention that the method 
of 
averaging in the first 
approximation
was utilized in [14,16] 
for studying
the asymptotic behavior of 
solutions of a particular class of systems of equations
with oscillatory decreasing 
coefficients.
\par
In the present paper we adopt 
the method 
of Shtokalo
for the problem of asymptotic 
integration of
systems of linear differential 
equations with oscillatory 
decreasing coefficients.
\par
We consider the following 
system of 
differential equations
in $n$-dimensional space 
${\cal R}^n$
$$
\frac{dx}{dt}=\left\{A_0+
\sum_{j=1}^k\frac{1}{t^{j\alpha}}
A_j(t)\right\} x+
\frac{1}{t^{(1+\delta)}}F(t)x.
\eqno(2)
$$
Here, $A_0$ is a constant 
$n\times n$ matrix, and
$A_1(t),A_2(t),..,A_k(t)$ are 
$n\times n$ matrices
that belong to $\Sigma$. 
We shall assume that matrix $A_0$ 
is in Jordan canonical form,
a real number $\alpha$ and 
a positive integer
$k$ satisfy 
$0<k\alpha\le 1<(k+1)
\alpha,~\delta>0$. The square
matrix $F(t)$ satisfies 
$$
||F(t)||\le C<\infty
$$
for $t_0\le t<\infty$, where 
$||.||$
is some matrix norm in 
${\cal R}^n$.
\par
We are concerned with the behavior 
of solutions
of system (2) when $t\to\infty$.
Let us construct an invertible
change of variables (for 
sufficiently large  
$t$,~$t>t^*>t_0$)
that would transform system 
(2) into  the simpler system
$$
\frac{dy}{dt}=\left\{\sum_{j=0}^k
\frac{1}{t^{j\alpha}}A_j\right\}y+
\frac{1}{t^{
(1+\varepsilon)}}G(t)y,
\quad\varepsilon>0,\quad t>t^*.
\eqno(3) 
$$
Here, 
$A_0,A_1,..,A_k$ are constants 
square matrices
(moreover, $A_0$ is the same 
matrix as in (1)),
and, the matrix $G(t)$ has the 
same 
properties as the matrix  
$F(t)$ in system (2).
\par
Without loss of generality, 
we can assume that
all eigenvalues of the matrix 
$A_0$ are real.
Indeed, if matrix $A_0$ has 
complex eigenvalues, 
then we can make a change of 
variables in (2)
$$
y=e^{iRt}z,
$$
where $R$ is diagonal matrix 
composed of the imaginary
parts of eigenvalues of the 
matrix
$A_0$.  This change of variable
with coefficients, which are 
bounded in
$t$,~$t\in (-\infty,\infty)$,
transforms the matrix 
$A_0$ into the matrix $A_0-iR$, which only has real eigenvalues.
\par
We shall try to choose an 
invertible
change of variables 
(for sufficiently large $t$), 
in the form
$$
x=\left\{\sum_{j=0}^k
\frac{1}{t^{j\alpha}}Y_j(t)
\right\}y, \eqno(4)
$$
to transform system (2) into 
system (3), 
where $Y_0(t)=I$ is the identity 
matrix, and, 
$Y_1(t),..,Y_k(t)$ are 
$n\times n$ matrices
that belong to $\Sigma$ and have 
zero mean value.
By substituting (4) into (2),
and, replacing
$\frac{dy}{dt}$ by the right 
hand side of (3) we obtain
$$
\begin{array}{l}
\left\{\sum_{j=0}^k\frac{1}
{t^{j\alpha}}Y_j(t)\right\}
\left\{\sum_{j=0}^k
\frac{1}{t^{j\alpha}}A_j\right\}
y+\\
+\frac{1}{t^{(1+\varepsilon)}}
\left\{\sum_{j=0}^k\frac{1}
{t^{j\alpha}}Y_j(t)
\right\}G(t)y+
\frac{1}{t^{(1+\alpha)}}W(t)y+
\left\{\sum_{j=1}^k\frac{1}
{t^{j\alpha}}\frac{Y_j(t)}{dt}
\right\}y=\\
=\left\{A_0+\sum_{j=1}^k
\frac{1}{t^{j\alpha}}A_j(t)
\right\}\left\{
\sum_{j=0}^k\frac{1}{t^{j\alpha}}
Y_j(t)\right\}y+
\frac{1}{t^{(1+\delta)}}U(t)y,
\end{array}\eqno(5)
$$
where
$$
W(t)=\left\{-\sum_{j=1}^k\frac
{j\alpha}{t^{(j-1)\alpha}}Y_j(t)
\right\},\eqno(6)
$$
and
$$
U(t)=F(t)\left\{\sum_{j=0}^k
\frac{1}{t^{j\alpha}}Y_j(t)
\right\}. \eqno(7)
$$
Equating the summands that 
contain
$t^{-j\alpha}$~$(j=1,\dots,k)$ 
in the left and the right hand 
sides of (5) yields a system
of $k$ linear matrix 
differential equation with 
constant coefficients
$$
\begin{array}{l}
\frac{dY_j(t)}{dt}-A_0Y_j(t)+
Y_j(t)A_0=\\
=\sum_{l=0}^{j-1}A_{j-l}(t)Y_l(t)-
\sum_{l=0}^{j-1}Y_l(t)A_{j-l},
\quad (j={1,\dots,k}).
\end{array}\eqno(8) 
$$
\par
The solvability of system (8)
was studied [11].  We 
represent $Y_j(t)$ as a finite 
sum
$$
Y_j(t)=\sum_{\lambda\ne 0}y_
{\lambda}^je^{i\lambda t},
$$
where $y_{\lambda}^j$ are 
constants $n\times n$
matrices, and, obtain matrix 
equations
$$
i{\lambda}y_{\lambda}^j-A_0y_
{\lambda}^j+y_{\lambda}^jA_0=
b_{\lambda}^j.
$$
Since all the eigenvalues of $A$ are real,
the matrix equations have 
unique solutions for $\lambda\ne 0$ 
(see, for instance, [16,17].
On each of the $k$ steps of 
the solution process
we determine the matrix $A_j$ 
from the condition that the 
right hand side
of (8) has a zero mean value. 
In particular, for $j=1$
$$
\frac{dY_1(t)}{dt}-A_0Y_1(t)+
Y_1(t)A_0=A_1(t)-A_1,
$$
where $A_1$ is the mean value of 
the matrix $A_1(t)$. 
Relation (5) implies the 
following result.

\bigskip

{\bf Theorem 1}. {\it System
(2), for sufficiently large $t$,
can be transformed using a change of variables
(4) into a system
$$
\frac{dy}{dt}=\left\{\sum_{j=0}^k
\frac{1}{t^{j\alpha}}A_j\right\}y
+\frac{1}{t^{(1+\varepsilon)}}
G(t)y,
$$
where $\varepsilon>0$, and 
$||G(t)||\le C_1<\infty$ . }

\bigskip

{\bf Proof}. Substituting (4) 
into (2), we obtain
$$
\begin{array}{l}\left\{
\sum_{j=0}^k\frac{1}
{t^{j\alpha}}Y_j(t)\right\}\frac
{dy}{dt}=\left\{A_0+
\sum_{j=1}^k\frac{1}
{t^{j\alpha}}
A_j(t)\right\}\left\{
\sum_{j=0}^k\frac{1}{t^{j\alpha}}
Y_j(t)\right\}y-\\
-\left\{\sum_{j=1}^k\frac{1}{t^
{j\alpha}}\frac{dY_j(t)}{dt}
\right\}y-
\frac{1}{t^{(1+\alpha)}}
W(t)y+
\frac{1}{t^{(1+\delta)}}U(t)y,
\end{array}
$$
where $W(t)$ and $U(t)$ are 
defined 
by (6) and (7) respectively. 
The last relation can be 
rewritten as
$$
\begin{array}{l}\left\{
\sum_{j=0}^k\frac{1}{t^
{j\alpha}}Y_j(t)\right\}
\left\{\frac{dy}{dt} -\left\{
\sum_{j=0}^k
\frac{1}{t^{j\alpha}}A_j\right\}y
\right\}=\\
=\left\{A_0+\sum_{j=1}^k\frac{1}
{t^{j\alpha}}A_j(t)\right\}
\left\{\sum_{j=0}^k\frac{1}
{t^{j\alpha}}Y_j(t)\right\}y-
\left\{\sum_{j=1}^k\frac{1}
{t^{j\alpha}}\frac{dY_j(t)}{dt}
\right\}y-\\
-\left\{\sum_{j=0}^k\frac{1}{t^
{j\alpha}}Y_j(t)\right\}
\left\{\sum_{j=0}^k\frac{1}{t^
{j\alpha}}A_j\right\}y-\frac{1}
{t^{(1+\alpha)}}
W(t)y+\frac{1}{t^{(1+\delta)}}
U(t)y. \end{array}
$$
Due to (8) we get
$$
\begin{array}{l}\left\{
\sum_{j=0}^k\frac{1}
{t^{j\alpha}}Y_j(t)\right\}
\left\{\frac{dy}{dt}-\left\{
\sum_{j=0}^k
\frac{1}{t^{j\alpha}}A_j\right\}
y\right\}=\\
=\frac{1}{t^{(1+k)\alpha}}S(t)y-
\frac{1}{t^{(1+\alpha)}}W(t)y+
\frac{1}{t^{(1+\delta)}}U(t)y,
\end{array}\eqno(9)
$$
where elements of the matrix
$S(t)$ can be represented as 
$t^{-j\alpha}a_j(t)$
$(j=0,\dots,k)$,
and $a_j(t)$ are trigonometric
polynomials.  Therefore,
$$
\frac{1}{t^{(1+k)\alpha}}S(t)-
\frac{1}{t^{(1+\alpha)}}W(t)+
\frac{1}{t^{(1+\delta)}}U(t)=
\frac{1}{t^{(1+\varepsilon)}}R(t),
\eqno(10)
$$
where $\varepsilon>0$ and 
$R(t)$ satisfies
$$
||R(t)||\le C_2<\infty.
$$
The identity (9), for 
sufficiently large $t$, can be 
rewritten
as
$$
\frac{dy}{dt}=\left\{\sum_{j=0}^k
\frac{1}{t^{j\alpha}}A_j\right\}y
+
\frac{1}{t^{(1+\varepsilon)}}
\left\{\sum_{j=0}^k\frac{1}
{t^{j\alpha}}
Y_j(t)\right\}^{-1}R(t)y.
$$
This relation and (10) yield
the assertion of the theorem.

\bigskip

The main part of system (3)
$$
\frac{dy}{dt}=\left\{\sum_{j=0}^k
\frac{1}{t^{j\alpha}}A_j\right\}y
$$
does not have oscillating 
coefficients.  This makes it
simpler than the original 
system (2).
In particular, the Fundamental 
Theorem of Levinson on asymptotic
behavior of solutions of 
linear systems of differential 
equations (see [1.2], as well as
[8,9]), readily yields the following
resolt for system (3).
%can be used for constructing the 
%asymptotics of the fundamental 
%matrix of the system (3). 
%We state Levinson's Theorem in a 
%form that is convenient
%for our consideration.

\bigskip

{\bf Theorem 2}. {\it Let us
first nonzero matrix among the
matrices $A_j$,~$j=0,1,\dots,k$
be the matrix $A_l$.
Let the matrix $A_l$ have distinct
eigenvalues.  Then the fundamental matrix
of system (3) has the following 
form
$$
X(t)=(P+o(1))exp\int_{t^*}^t 
\Lambda(s)ds,\quad t>t^*,\quad
t\to \infty,
$$
where $P$ is a matrix composed 
of the eigenvectors of the matrix
$A_l$, and $\Lambda(t)$ is a 
diagonal matrix
whose elements are eigenvalues 
of the matrix
$\sum_{j=l}^k t^{-j\alpha}A_j$.}

\bigskip

To prove this theorem we just 
have to observe
that the system of differential 
equations
$$
\frac{dx}{dt}=\frac{1}{t^l}A_lx+
\sum_{j=l+1}\frac{1}{t^j}A_jx
$$
can be transformed into 
$$
\frac{dx}{d\tau}=\frac{1}{1-l}[A_l+
\sum_{j=l+1}\frac{1}{t^{j-l}}A_j]x
$$
using the change of variables 
$\tau=t^{1-l}$.

\bigskip

As an example we consider an 
equation of an adiabatic oscillator
$$
\frac{d^2y}{dt^2}+(1+\frac{1}
{t^{\alpha}}\sin \lambda t)y=0,
 \eqno(11)
$$
where $\lambda$,~$\alpha$ are 
real numbers,
and $0<\alpha\le 1$.
The problem of asymptotic integration of
the equation (11) has been 
studied in [5,6, 18--20]. 
In particular, asymptotics of 
solutions for 
$\frac{1}{2} \le \alpha \le 1$ 
were obtained. 
The method that we proposed in 
this chapter
can be used to obtain 
(in a simple manner)
all known results on asymptotics 
of solutions
of equation (11) as well as to 
establish new results.
\par
Let us pass from equation (11) 
to the 
system of
equations $(x=(x_1,x_2))$ using a 
change of variables
$$
y=x_1 \cos t + x_2 \sin t, \quad
y'= -x_1 \sin t + x_2 \cos t.      
\eqno(12)
$$
We obtain the system
$$
\frac{dx}{dt}=\frac{1}
{t^{\alpha}}A(t)x.
\eqno(13)
$$
It is convenient to rewrite
the matrix $A(t)$ in complex 
form as
$$
\aligned
&A(t)=a_1e^{i(\lambda+2)t}+
\bar a_1e^{-i(\lambda+2)t}+
a_2e^{i(\lambda-2)t}+
\bar a_2e^{-i(\lambda-2)t}+ 
a_3e^{i\lambda t}+\bar a_3
e^{-i\lambda t},
\endaligned
$$
where
$$
a_1=\frac{1}{8}
\left(\begin{array}{cc}
-1 & i\\
i & 1
\end{array}\right),\quad
a_2=\frac{1}{8}
\left(\begin{array}{cc}
1 & i\\
i &-1
\end{array}\right),\quad
a_3=\frac{1}{8}
\left(\begin{array}{cc}0&-2i
\\2i&0
\end{array}\right),
$$
and the matrices $\bar a_1$,~
$\bar a_2$,~$\bar a_3$ are complex
conjugates to the matrices 
$a_1,a_2,a_3$ respectively.
\par
The values of $\alpha$ and 
$\lambda$
significantly affect the behavior of solutions
of system (13). We denote by 
$R(t)$ a $2\times 2$ matrix that satisfies
$$
||R(t)||\le C_3 <\infty.
$$
for all $t$.
\par
First, assume
$\frac{1}{2}<\alpha\le 1$. 
For $\lambda\ne\pm 2$
system (3) becomes
$$
\frac{dy}{dt}=\frac{1}
{t^{1+\varepsilon}}R(t)y,\quad
\varepsilon>0.
$$
Therefore, it is easy to see 
(taking into consideration 
change (12)),
that the fundamental system of
solutions
of equation (11) for
$\frac{1}{2}<\alpha\le 1$,~ 
$\lambda\ne\pm 2$ as $t\to\infty$
has the form
$$
x_1= \cos t +o(1),\quad x_2=
\sin t +o(1),
$$
$$
x'_1= -\sin t +o(1),\quad x'_2=
\cos t + o(1).
$$
We shall represent the fundamental system
of solutions of equation (11) 
as a matrix with rows
$x_1,x_2$ and $ x'_1,x'_2$.
\par
Now assume $\lambda=\pm2$.  
More 
specifically let $\lambda=2$.
Then system (3) becomes
$$
\aligned
&\frac{dy}{dt}=\frac{1}
{t^{\alpha}}A_1y+\frac{1}
{t^{1+\varepsilon}}R(t)y,
\quad\varepsilon>0.
\endaligned 
$$
Here
$$
a_2 +\bar a_2=
A_1=
\frac{1}{4}\left
(\begin{array}{cc}
1 & 0\\
0 & -1
\end{array}\right).
$$
\par
Theorem 2 implies that for 
$t\to\infty$ the fundamental matrix
of system (11) has the following
form
$$
Y(t)=\left(
\begin{array}{cc}
exp\left(\int_{t^*}^{t}
\frac{1}{4s^{\alpha}}ds
\right) & 0 \\
0 & exp\left(-\int_{t^*}^
{t}\frac{1}{4s^{\alpha}}ds
\right)
\end{array}\right)
\left[I+o(1)\right].
$$
Therefore, for $\alpha=1,~\lambda=2$ we obtain
the fundamental system of solutions 
of equation (11) as 
$$
\begin{array}{l}\left
(\begin{array}{cc}
t^{\frac{1}{4}}\cos t & 
t^{-\frac{1}{4}}\sin t\\
-t^{\frac{1}{4}}\sin t & 
t^{-\frac{1}{4}}\cos t
\end{array}\right)
\left[I+o(1)\right],
\end{array}
$$
while for $1/2<\alpha<1,
~\lambda=2$ 
for $t\to\infty$ we get
$$ 
\begin{array}{l}\left(\begin{array}{cc}
exp\left(\frac{t^{1-\alpha}}
{4(1-\alpha)}\right)\cos t &
exp\left(-\frac{t^{1-\alpha}}
{4(1-\alpha)}\right)\sin t \\
-exp\left(\frac{t^{1-\alpha}}
{4(1-\alpha)}\right)\sin t &
exp\left(-\frac{t^{1-\alpha}}
{4(1-\alpha)}\right)\cos t
\end{array}\right)
\left[I+o(1)\right].
\end{array}
$$
\par
We note that, for $\lambda=\pm 2$, 
and $\frac{1}{2}<\alpha\le 1$,  
equation (11)
has unbounded solutions. Moreover,
for $\alpha=1$, the solutions 
have a polynomial growth,
while, for $\alpha\ne 1$, they 
grow exponentially.
\par
We now assume 
$\frac{1}{3}<\alpha
\le\frac{1}{2}$.
In this case a change of variables
(4) transforms (13) into
$$
\frac{dy}{dt}=\frac{1}
{t^{\alpha}}A_1y+\frac{1}
{t^{2\alpha}}A_2y+\frac{1}
{t^{1+\varepsilon}}R(t)y,
\quad\varepsilon>0.
$$
If $\lambda\ne\pm 2,\pm 1$, then 
the matrix $A_1$ is zero, and
matrix $A_2$ has the form
$$
A_2 = i\left[\frac{1}{{\lambda + 2}}
(a_1\bar a_1 - \bar a_1
 a_1)+
\frac{1}{{\lambda -2}}(a_2
\bar a_2 -\bar a_2 a_2)
 +\frac{1}{\lambda}(a_3
\bar a_3 -\bar a_3 
a_3)\right].  \eqno(14)
$$
Computing $A_2$ yields
$$
A_2=\frac{1}{4(\lambda^2-4)}
\left(\begin{array}{cc}
0 & 1\\
-1 & 0
\end{array}\right).
$$
The system 
$$
\frac{dy}{dt}=\frac{1}{t^
{2\alpha}}A_2y
$$
can be integrated.  We obtain
that, for $t\to\infty$, the
fundamental system of solutions
of equation (11) with
$\alpha=\frac{1}{2}$,
~$\lambda\ne\pm 2,\pm 1$,
has the form
$$
\left(\begin{array}{cc}
\cos(t+\gamma\ln t) &\sin(t+
\gamma\ln t)\\
-\sin(t+\gamma\ln t) & \cos
(t+\gamma\ln t)
\end{array}\right)[I+o(1)],
$$
where $\gamma=\frac{1}{4
(\lambda^2-4)}$.  For
$\frac{1}{3}<\alpha<\frac{1}{2}$,
and $\lambda\ne\pm 2,\pm 1$, 
the fundamental system of solutions
of equation (11) has the form
$$
\left(\begin{array}{cc}
\cos\left(\frac{t^{1-2\alpha}}
{4(1-2\alpha)(\lambda^2-4)}
\right) &
\sin\left(\frac{t^{1-2\alpha}}
{4(1-2\alpha)(\lambda^2-4)}
\right)\\
-\sin\left(\frac{t^{1-2\alpha}}
{4(1-2\alpha)(\lambda^2-4)}
\right) &
\cos\left(\frac{t^{1-2\alpha}}
{4(1-2\alpha)(\lambda^2-4)}
\right)
\end{array}\right)[I+o(1)]
$$
as $t\to\infty$.
\par
We now assume 
$\alpha=\frac{1}{2},
\lambda=1$.  In this case
$A_1$ is zero, and $A_2$ is determined by
$$
iA_2=-\frac{1}{3}a_1\bar a_1
+\frac{1}{3}\bar a_1a_1-
\bar a_2a_2+a_2\bar a_2
-a_3\bar a_3 +\bar a_3 
a_3+
a_2a_3+a_3a_2-\bar a_2
\bar a_3+\bar a_3
\bar a_2 .
$$
A simple calculation yields
$$
A_2=\frac{1}{24}
\left(\begin{array}{cc}
0 & -5\\
-1 & 0
\end{array}\right).
$$
The corresponding system (3)
has the following form
$$
\frac{dy}{dt}=\frac{1}
{t}A_2y+\frac{1}{t^{1+
\varepsilon}}R(t)y,\quad
\varepsilon>0.
$$
By integrating the system
$$
\frac{dy}{dt}=\frac{1}
{t}A_2y,   
$$
we obtain its fundamental matrix
$$
Y(t)=\left(\begin{array}{cc}
-\sqrt 5 t^\varrho & \sqrt 5 
t^{-\varrho}\\
t^\varrho & t^{-\varrho}
\end{array}\right),
$$
where $\varrho=\frac{\sqrt 5}
{24}$. Then the fundamental
system of solutions of equation
(11) for
 $\alpha=\frac{1}{2}$, 
$\lambda=1$, and $t\to\infty$,
has the form
$$
\left(\begin{array}{cc}
t^{\varrho}\sin(t-\beta) & t^{-\varrho}\sin(t+\beta)\cr
t^{\varrho}\cos(t-\beta) & 
t^{-\varrho}\cos(t+\beta)
\end{array}\right)[I+o(1)],
$$
where
$$
\varrho=\frac{\sqrt 5}{24},
\quad  \beta=arctg\sqrt{5},\quad 
0<\beta<\frac{\pi}{2}. 
\eqno(15)
$$
If $\frac{1}{3}<\alpha<\frac{1}{2}$ and $\lambda=1$ 
we have the system 
$$
\frac{dy}{dt}=\frac{1}
{t^{2\alpha}}A_2y+\frac{1}
{t^{1+\varepsilon}}R(t)y,
\quad\varepsilon>0.
$$
Using Theorem 2 we obtain
the asymptotics of the fundamental matrix of this system,
and then, using the change 
(12), the asymptotics of
the fundamental system of solutions of equation (1.7.12) for
$\frac{1}{3}<\alpha<\frac{1}{2}$, $\lambda=1$,
and $t\to\infty$:
$$
\left(\begin{array}{cc}
exp(\varrho\frac{t^{1-2\alpha}}
{1-2\alpha})
\sin (t-\beta) &
exp(-\varrho\frac{t^{1-2\alpha}}
{1-2\alpha})\sin (t+\beta)
\cr
exp(\varrho\frac{t^{1-2\alpha}}
{1-2\alpha})\cos (t-\beta)
& exp(-\varrho\frac{t^{1-2\alpha}}
{1-2\alpha})\cos (t+\beta)
\end{array}\right)[I+o(1)],
$$
where $\varrho$ are $\beta$  
are defined by (15).
Thus, for $\alpha=\frac{1}{2}$ 
and $\lambda=1$
we observe a polynomial growth 
of solutions, while for 
$\frac{1}{3}<\alpha<\frac{1}{2}$ 
and $\lambda=1$ the
solutions grow exponentially.
\par
Now let $\alpha=\frac{1}{2}$ 
and $\lambda=2$. Simple calculations
show that 
$$
A_1=\left(\begin{array}{cc}
\frac{1}{4} & 0\\
0 & -\frac{1}{4}
\end{array}\right),\quad
A_2=\left(\begin{array}{cc}
0 & -\frac{1}{64}\\
\frac{1}{64} & 0
\end{array}\right).
$$
Therefore, we get a system
$$
\frac{dy}{dt}=\frac{1}{t^{\frac{1}
{2}}}A_1y+\frac{1}{t}A_2y+
\frac{1}{t^{1+
\varepsilon}}R(t)y,\quad
\varepsilon>0. \eqno(16)
$$
We compute the eigenvalues of 
the matrix
$$
\frac{1}{t^{\frac{1}{2}}} A_1
+\frac{1}{t}A_2,
$$ 
integrate them, and, using 
Theorem 2,
we obtain the asymptotics of the
fundamental matrix of system 
(16).
Next we find the fundamental system of
solutions of equation (11) for $\alpha=\frac{1}{2},\,
\lambda=2$ and $t\to\infty$:
$$
\left(\begin{array}{cc}
exp(\phi(t))\cos t & 
exp-(\phi(t))\sin t\\
-exp(\phi(t))\sin t & 
exp-(\phi(t))\cos t
\end{array}\right)[I+o(1)],
$$
where $\phi(t)=\frac{1}{2}
\sqrt{t}$.  
%The last formulas
%give more precise asymptotics than
%the corresponding formulas in
%Harris, Lutz [1977].
\par
For $\frac{1}{3}<\alpha<
\frac{1}{2},~\lambda=2$ 
instead of (16) we get a system
$$
\frac{dy}{dt}=\frac{1}{t^{\alpha}}
A_1y+\frac{1}{t^{2\alpha}}A_2y+
\frac{1}{t^{1+\varepsilon}}
R(t)y,
$$
where $\varepsilon>0$, with the same matrices $A_1$, $A_2$. 
Therefore, it is straightforward to write the asymptotics
of the fundamental system of solutions
of equation (11) for 
$\frac{1}{3}<\alpha<\frac{1}{2}$, $\lambda=2$, and $t\to\infty$.
\par
Finally, let 
$\alpha=\frac{1}{3}$,~
$\lambda\ne\pm 1$, and 
$\lambda\ne\pm 2$.
Then, it turns out that 
$A_1$ is zero, and $A_2$ is 
defined by (14).
Matrix $A_3$ differs from 
zero only if
$\lambda=\pm\frac{2}{3}$. 
Assume $\lambda=\frac{2}{3}$.
System (3) then becomes
$$
\frac{dy}{dt}=\frac{1}{t^{\frac{2}
{3}}}A_2 y+\frac{1}{t}A_3y + 
\frac{1}{t^
{1+\varepsilon}}R(t)y,
$$
where $\varepsilon>0$ and the matrices $A_2$ and $A_3$ 
are defined by
$$
A_2=\left(\begin{array}{cc}
0 & -\frac{9}{128}\\
\frac{9}{128} & 0
\end{array}\right),\quad
A_3=\left(\begin{array}{cc}
-\frac{27}{1024} & 0\\
0 & \frac{27}{1024}
\end{array}\right).
$$
We compute the eigenvalues of 
the matrix
$$
\frac{1}{t^{\frac{2}{3}}}A_2+
\frac{1}{t}A_3. 
$$
These eigenvalues have zero real parts,
for sufficiently large $t$.
Further, using the same 
scheme as before
we find the asymptotics of 
the fundamental
system of solutions of 
equation (11).
We only note that the solutions of equation
(11) are bounded for
$\alpha=\frac{1}{3},\,
\lambda=\pm\frac{2}{3}$ as
$t\to\infty$.

\bigskip

\begin{center}{\bf References}
\end{center}

1. N. Levinson  "The asymptotic nature of solution of linear differential
equations," Duke Math. J.,  
15,111 -- 126, (1948).

2. I.M. Rappoport, On some 
asymptotic methods in the theory
of differential equations [in Russian],
Izd. AN UssR, Kiev, (1954).

3. P. Hartman and A. Wintner 
"Asymptotic integration of 
linear differential
equation," Amer. J. Math.,  
77, 45 -- 86, 932, (1955).

4. M.V. Fedoryuk   
"Asymptotic methods in theory
of one-dimensional singular 
differential operators,"
Trudy Moscov. Mat. Obshch. 
[Transl. Moskow. Math. Soc.], 
15, 296 -- 345, (1966).

5. W.A. Harris, Jr. and D.A.
Lutz  "On the asymptotic integration of
linear differential systems,"  
 J. Math. Anal. and Applic.,
 1974. 48. No. 1,1 -- 16, (1974).

6. W.A. Harris, Jr. and D.A. Lutz  
"A Unified theory of asymptotic 
integration,"
 J. Math. Anal. and Applic., 
57, 571 -- 586, (1971)

7. R. Bellman, Stability Theory
of Differential Equations, New York
(1953).

8. E. Coddington and N. Levinson,
Theory of Ordinary Differential
Equations, New York (1955).

9. M.A. Naimark, Linear Differential
Operators [in Rossian], Nauka, Moscow
(1969).

10. L. Cezari, Asymptotic Behavior
and Stability Problems
of Ordinary Differential 
Equations, West Berlin (1959).

11. M.S.P. Eastham, The 
asymptotic solution of linear 
differential systems,
London Math. Soc. Monographs,
Vol 4. Clarendon Press, (1989).

12. I.Z.Shtocalo, "Tect for stability
and instability of the solutions
of linear differential equations
with almost periodic coefficients,"
Math, Sb. [Math. USSR Sb.], 
19(61), No.2,263 -- 286 (1946).

13. I.Z. Shtokalo, Linear 
Differential Equations With
Variable Coefficients [in Russian],
Izd. AN USSR, Kiev(1960).

14. Yu.A. Samokhin and V.N. Fomin,
"A method for studying the stability
of the solutions of linear of system
subjected to the action of parametric
loads with continuous spectrom,"
Sibirsk. Mat. Zh. [Siberian Math.
J.], 17, No.4., 926 -- 931(1976).

15. Yu.A. Samokhin and V.N. Fomin,
The asymptotic integration of system
of differential equations with
oscillating decreasing coefficients,"
in: Problems of the Theory 
Of Periodic Motions [in Russian],
Vol. 5, Izhevsk(1981), pp. 45--50.

16. F.R. GantmakherГан, The Theory
of Matrices [in Russian], Nauka,
Moscow(1966).

17. Yu. A. Daletskii and M.G.
Krein, Stability of the 
Solution of Differential Equations
in Banach Space [in Russian], Nauka,
Moscow(1970).

18. A. Wintner "The adiabatic 
linear oscillator,"  Amer. J. 
Math., 68, 385 -- 397,
(1946).

19. A.Wintner "Asymptotic 
integration of the  adiabatic 
oscillator," Amer. J. Math.,
69. 251 -- 272, (1946).

20. W.A. Harris,Jr. and D.A.
Lutz "Asymptotic integration  
of adiabatic oscillator,"
J. Math. Anal. and Applic.,
51, No. 1, 76 -- 93, (1975).

\end{document}